\documentclass[12pt]{amsart}
\usepackage{amsfonts,amssymb, fullpage, color,enumerate, url}

\usepackage[utf8]{inputenc}   
\usepackage[T1]{fontenc}

\def\Pr{\begin{proof}}
\def\Rp{\end{proof}}

\def\norm{\left\Vert{.}\right\Vert}
\newcommand\norme[1]{\|#1\|}
\def\IR{\mathbb R}
\def\IN{\mathbb N}
\def\IQ{\mathbb Q}
\def\IK{\mathbb K}

\def\IL{\mathbb L}
\def\IK{\mathbb K}
\def\IF{\mathbb F}

\def\ZF{\mathbf {ZF}}
\def\ZFA{\mathbf {ZFA}}
\def\HB{\mathbf {HB}}
\def\BPI{\mathbf {BPI}}
 \def\A{\mathbf A}
\def\AC{\mathbf {AC}}

\def\B{\mathbf {B}}
\def\D{\mathbf {D}}

\def\MC{\mathbf {MC}}

\def\BE{\mathbf {BE}}
\def\LE{\mathbf {LE}}

\def\I{\mathbf {I}} 

\DeclareMathOperator{\spanv}{span}

\DeclareMathOperator{\Id}{Id}
\DeclareMathOperator{\BL}{BL}
\DeclareMathOperator{\can}{can}

\theoremstyle{plain}
\newtheorem{theorem}{Theorem}[section] 
\newtheorem{corollary}[theorem]{Corollary}
\newtheorem{proposition}[theorem]{Proposition}
\newtheorem{lemma}[theorem]{Lemma}
\newtheorem{remark}[theorem]{Remark}

\newtheorem*{proposition*}{Proposition}
\newtheorem*{theorem*}{Theorem}
\newtheorem*{corollary*}{Corollary}
\newtheorem*{lemma*}{Lemma}

\theoremstyle{definition}
\newtheorem{definition}[theorem]{Definition}
\newtheorem{notation}[theorem]{Notation}
\newtheorem{question}[theorem]{Question}

\theoremstyle{remark}
\newtheorem{example}[theorem]{Example}

\date{\today}
\usepackage{graphicx}

\begin{document}
\title{Linear extenders and the Axiom of Choice}
\author[M.~Morillon]{Marianne Morillon}
 \address{Laboratoire d'Informatique et Mathématiques, 
Parc Technologique Universitaire, Bâtiment 2,
2 rue Joseph Wetzell, 97490 Sainte Clotilde}
 \email{Marianne.Morillon@univ-reunion.fr}
 \urladdr{http://lim.univ-reunion.fr/staff/mar/}

 \begin{abstract}  In set theory without the axiom of Choice $\ZF$, we prove that for every 
 commutative field $\IK$, the following statement $\D_{\IK}$: ``On every non null $\IK$-vector 
space, there exists a non null linear form''  implies the existence of a ``$\IK$-linear extender'' 
on every vector subspace of a $\IK$-vector space. This solves a question raised in \cite{Mo09}.
 In the second part of the paper, we generalize our results in the case of spherically complete 
 ultrametric valued fields, and show that Ingleton's statement is equivalent to the existence of 
``isometric linear extenders''. 
 \end{abstract}

 \maketitle

\section{Introduction}
We work in $\ZF$, set theory without the Axiom of Choice ($\AC$). 
Given a commutative field $\IK$ and two  $\IK$-vector spaces $E$, $F$, 
we denote by $L_{\IK}(E,F)$ (or $L(E,F)$) the set of $\IK$-linear mappings $T: E \to F$. Thus, 
$L(E,F)$ is a vector subspace of the product vector space  $F^E$. 
 A {\em linear form} on the $\IK$-vector space $E$
is a $\IK$-linear mapping $f: E \to \IK$; we denote by $E^*$ the {\em algebraic dual} of $E$,
{\em i.e.} the vector space $L(E,\IK)$. 
Given two  $\IK$-vector spaces $E_1$ and $E_2$, and a linear mapping $T: E_1 \to E_2$, 
we denote by $T^t: E_2^* \to E_1^*$ the mapping associating to every $g \in E_2^*$ the linear 
mapping $g \circ T: E_1 \to \IK$: the mapping $T^t$ is $\IK$-linear and is called the {\em 
transposed} mapping of $T$. 
Given a vector subspace $F$ of the $\IK$-vector space $E$, a 
{\em linear extender} for $F$ in $E$  is a linear mapping $T: F^* \to E^*$ associating to each linear form
$f \in F^*$ a linear form $\tilde f: E \to \IK$ extending $f$. If 
  $G$ is a {\em complementary subspace} of $F$ in $E$ ({\em i.e.} $F + G= E$ and $F \cap G=\{0\}$), 
  which we denote by  $F \oplus G= E$, if $p: E \to F$ is the linear mapping fixing every element 
of $F$ and which is null on $G$, then the transposed mapping  $p^t: F^* \to E^*$ associating
to every $f \in F^*$ the linear form $f \circ p: E \to \IK$ is a linear extender for $F$ in $E$. 
However, given a commutative field $\IK$, the existence of a complementary subspace of every 
subspace of a $\IK$-vector space implies $\AC$ in $\ZF$. More precisely, 
denoting  by $\ZFA$ (see \cite[p.~44]{Je73}) the set theory $\ZF$ with the axiom of  
extensionality 
weakened to allow the existence of atoms, 
 the existence of a complementary subspace for every subspace of a $\IK$-vector space 
  implies,  in $\ZFA$ 
 (see \cite[Lemma~2]{Blei}), the following  {\em Multiple Choice} axiom $\MC$ (see 
\cite[p.~133]{Je73} and  form~37 of \cite[p.~35]{Ho-Ru}): {\em ``For every infinite family $(X_i)_{i 
\in I}$ of nonempty sets, there exists a family 
$(F_i)_{i \in I}$ of nonempty finite sets such that for each $i \in I$, $F_i \subseteq X_i$.''}
It is known that $\MC$ is equivalent to $\AC$ in $\ZF$, but $\MC$ does not imply $\AC$ in 
$\ZFA$. 

 Given a commutative field $\IK$, we consider the following consequences of the Axiom of Choice:
\begin{itemize}
\item  $\BE_{\IK}$: {\em ``Every linearly independent subset of a vector space $E$ over $\IK$ is 
included in  a basis of $E$.''}
\item  $\B_{\IK}$: {\em ``Every vector space over $\IK$ has a basis.''}
\item  $\LE_{\IK}$: (Linear Extender) {\em ``For every subspace $F$ of a $\IK$-vector space $E$, 
there 
exists a linear mapping $T:F^* \to E^*$ associating to every $f \in F^*$ a linear mapping $T(f): E 
\to \IK$ extending $f$.}
\item  $\D_{\IK}$: {\em ``For every non null $\IK$-vector space $E$ there exists a non null linear 
form $f:E \to \IK$.}
\end{itemize}
 In $\ZF$,   $\BE_{\IK} \Rightarrow \B_{\IK} \Rightarrow \LE_{\IK} \Rightarrow \D_{\IK} $ 
(see \cite[Proposition 4]{Mo09}). 
In this paper, we  show (see Theorem~\ref{theo:d_to_le} in Section~\ref{sec:d_to_le}) 
that for each commutative field $\IK$, $\D_{\IK}$ implies $\LE_{\IK}$, and this solves 
Question~2 of \cite{Mo09}. 
In Section~\ref{sec:other_equ} we provide several other statements which are equivalent to 
$\D_{\IK}$ and we  
introduce the following consequence $w\D_{\IK}$ of $\D_{\IK}$: 
 {\em ``For every  $\IK$-vector space $E$ and every 
$a \in E \backslash \{0\}$, there exists an additive mapping $f: E \to \IK$ 
such that  $f(a)=1$.''} 

\begin{question} \label{quest:add_to_lin}
Given a commutative field $\IK$, does  the statement $w\D_{\IK}$  imply $\D_{\IK}$?
\end{question}

\medskip

Blass (\cite{Bla84}) has shown that the statement $\forall \IK  \B_{\IK}$ (form~66 of \cite{Ho-Ru}: 
``For every commutative field,  every $\IK$-vector space has a basis'')    implies $\MC$ in 
$\ZFA$ (and thus  implies $\AC$ in $\ZF$), but it is an open question to know whether there exists 
a commutative field $\IK$ such that $\B_{\IK}$ implies $\AC$. 
In $\ZFA$, the statement $\MC$ implies $\D_{\IK}$ for every commutative field $\IK$ with null 
characteristic 
(see \cite[Proposition~1]{Mo09}). Thus in $\ZFA$, the statement ``For every commutative field $\IK$ 
with null characteristic, $\D_{\IK}$''   does not imply $\AC$.
Denoting by $\BPI$  the {\em Boolean prime ideal}: ``Every non null  boolean algebra has an 
ultrafilter'' 
(see form~14 in \cite{Ho-Ru}), Howard and Tachtsis (see \cite[Theorem~3.14]{How-Tach13}) have shown 
that  for every finite field $\IK$,   $\BPI$ implies $\D_{\IK}$. Since $\BPI \not \Rightarrow 
\AC$, the statement ``For every  finite field $\IK$, $\D_{\IK}$'' does  not imply $\AC$. 
They also have shown  (see \cite[Corollary~4.9]{How-Tach13})  that  in $\ZFA$, $\forall \IK  
\D_{\IK}$ 
(``For every commutative field, for every non null $\IK$-vector space $E$, there  exists a 
$\IK$-linear form $f: E \to \IK$'') 
does not imply $\forall \IK  \B_{\IK}$, however, the following questions  seem to be  open in 
$\ZF$:
\begin{question}
Does the statement $\forall \IK \D_{\IK}$ imply $\AC$ in $\ZF$?
Is there  a (necessarily infinite) commutative field $\IK$ 
such that  $\D_{\IK}$ implies  $\AC$ in $\ZF$? 
\end{question}
 
In Section~\ref{sec:ingleton}, we extend Proposition~1 of  Section~\ref{sec:d_to_le} to the case of 
  spherically complete ultrametric valued fields (see Lemma~\ref{lem:equiv-ing}) and prove that 
Ingleton's statement, which is a ``Hahn-Banach type'' result for 
ultrametric semi-normed spaces over spherically complete ultrametric valued fields $\IK$, 
follows from $\MC$ when $\IK$ has a null characteristic. In Section~\ref{sec:cont-le}, we prove that 
Ingleton's statement is equivalent to the existence of ``isometric linear extenders''.

\section{$\LE_{\IK}$ and  $\D_{\IK}$ are equivalent} \label{sec:d_to_le}
\subsection{Reduced powers of a commutative field $\IK$}
Given a set $E$ and a filter $\mathcal F$ on a set $I$, 
we denote by $E_{\mathcal F}$ the quotient of the set $E^I$ by the equivalence relation $=_{\mathcal 
F}$
on $E^I$ satisfying  for every $x=(x_i)_{i \in I}$ and  $y=(y_i)_{i \in I} \in E^I$, 
$x =_{\mathcal F} y$ if and only if $\{i \in I : x_i=y_i\} \in \mathcal F$. If $\IL$ is a first 
order language and if $E$ carries a $\IL$-structure, then the quotient set $E_{\mathcal F}$ also 
carries a quotient $\IL$-structure: this $\IL$-structure   is a  {\em reduced power} of the 
$\IL$-structure $E$ (see \cite[Section~9.4]{Ho93}). Denoting by $\delta: E \to E^I$ the {\em 
diagonal mapping} associating to each
$x \in E$ the constant family $i \mapsto x$, and denoting by 
$\can_{\mathcal F}: E^I \to E_{\mathcal F}$ the canonical quotient mapping, then 
we denote by $j_{\mathcal F}: E \to E_{\mathcal F}$ the one-to-one mapping 
$\can_{\mathcal F} \circ \delta$. Notice that  $j_{\mathcal F}$ is a morphism of $\IL$-structures. 

 \begin{example}[The reduced power $\IK_{\mathcal F}$ of a field $\IK$]
Given a commutative field $\IK$ and a filter $\mathcal F$ on  a set $I$,  
we consider the  unitary $\IK$-algebra $\IK^I$, 
then the quotient $\IK$-algebra $\IK_{\mathcal F}$ 
is the quotient of the $\IK$-algebra  $\IK^I$ by the following ideal 
$nul_{\mathcal F}$   of {\em $\mathcal F$-almost everywhere null} elements of $\IK^I$:
$\{x=(x_i)_{i \in I} \in \IK^I : \{i \in I : x_i=0\} \in \mathcal F\}$. 
The mapping $j_{\mathcal F} : \IK \to \IK_{\mathcal F}$ is a one-to-one unitary morphism of 
$\IK$-algebras,  thus $\IK$ can be viewed as the one-dimensional unitary  $\IK$-subalgebra  of the 
$\IK$-algebra $\IK_{\mathcal F}$. Notice that the $\IK$-algebra $\IK_{\mathcal F}$ is a field if and 
only if $\mathcal F$ is an ultrafilter.
\end{example}

\begin{notation}
For every set $E$, we denote by $fin(E)$ the set of finite subsets of $E$, and we denote by 
$fin^*(E)$ the set of nonempty finite subsets of $E$.
\end{notation}

Given two sets $E$ and $I$, a binary relation  $R \subseteq E \times I$ is said to be  
{\em concurrent} (see \cite{Lu}) if for every $G \in fin^* (E)$, the set $R[G]:=\cap_{x \in G}R(x)$
is nonempty;
in this case, $\{R(x) :  x \in I\}$  satisfies the finite intersection property, and we denote by 
$\mathcal F_R$ the filter on $I$ generated by the sets $R(x)$, $x \in E$.

\subsection{$\D_{\IK}$ implies linear extenders}
\begin{remark} \label{rem:equ_DK}
It is known (see \cite[Theorem~2]{Mo09}), that $\D_{\IK}$ is equivalent to the following statement: 
``For every vector subspace of a $\IK$-vector space $E$ and every linear mapping
$f: F \to \IK$, there exists a $\IK$-linear mapping $f: E \to \IK$ extending $f$.''
\end{remark}

\begin{proposition} \label{prop:D_red_puis}
Given a commutative  field $\IK$, the following statements are equivalent:
\begin{enumerate}[i)]
\item $\D_{\IK}$
\item \label{it:dkbis} For every filter $\mathcal F$ on a set $I$, the linear  embedding 
$j_{\mathcal F} : \IK \to {\IK_{\mathcal F}}$ has a $\IK$-linear retraction 
$r:{\IK_{\mathcal F}} \to \IK$.
\end{enumerate}
\end{proposition}
\Pr $\Rightarrow$ Let $f: j_{\mathcal F}[\IK] \to \IK$ be the  mapping 
$x \mapsto j_{\mathcal F}^{-1}(x)$. Then $f$ is $\IK$-linear and $j_{\mathcal F}[\IK]$ is a vector 
subspace of the $\IK$-vector space ${\IK}_{\mathcal F}$. Using Remark~\ref{rem:equ_DK}, let $\tilde 
f: {\IK}_{\mathcal F} \to \IK$ be a $\IK$-linear mapping extending $f$; then $\tilde f$
is a $\IK$-linear retraction of $j_{\mathcal F}: \IK \to {\IK}_{\mathcal F}$. \\
$\Leftarrow$ Let $E$ be a non null $\IK$-vector space. Let $a$ be a non-null element of  $E$.
Using  Lemma~1 in \cite{Mo09},  there exists a filter $\mathcal F$ on the set $I=\IK^E$, and   a 
$\IK$-linear mapping 
$g : E \to \IK_{\mathcal F}$ such that $g(a)=j_{\mathcal F}(1)$.
Using~\eqref{it:dkbis}, let $r:  \IK_{\mathcal F} \to \IK$ be a $\IK$-linear retraction of 
the linear embedding $j_{\mathcal F}: \IK \to  \IK_{\mathcal F}$.
It follows that 
$f=r \circ g:E \to \IK$ is a $\IK$-linear mapping such that $f(a)=1$.
\Rp

\begin{theorem} \label{theo:d_to_le}
$\D_{\IK} \Leftrightarrow \LE_{\IK}$.
\end{theorem}
\Pr The implication $\LE_{\IK} \Rightarrow \D_{\IK}$ is trivial. We shall prove  $\D_{\IK} 
\Rightarrow \LE_{\IK}$. Given some vector subspace $F$ of a vector space $E$, let $I$ be the set of 
mappings 
$\Phi: F^* \to E^*$ and let $R$ be the binary relation on $fin(F^*) \times I$ such that 
for every $Z  \in fin(F^*)$ and every $\Phi \in I$, $R(Z,\Phi)$ if and only if for every
$f \in Z$, the linear form  $\Phi(f)$ extends $f$ and $\Phi$ is $\IK$-linear on 
$\spanv_{F^*}(Z)$.
Then the binary relation $R$ is concurrent: given $m$ finite subsets 
$Z_1, \dots, Z_m  \in fin(F^*)$, let $B=\{f_1,\dots,f_p\}$ be a (finite) basis of the $\IK$-vector 
subspace of $F^*$ generated by the finite set $\cup_{1 \le i \le m}Z_i$; 
then using $\D_{\IK}$ (see Remark~\ref{rem:equ_DK}), let 
$\tilde f_1$, \dots, $\tilde f_p$ be linear forms on $E$ extending $f_1$, \dots, $f_p$;
let $L: \spanv(\{f_1,\dots,f_p\}) \to E^*$ be the linear mapping such that for each 
$i \in \{1,\dots,p\}$, $L(f_i)=\tilde{f_i}$, and let 
$\Phi: F^* \to E^*$ be some mapping extending $L$ (for example, define $\Phi(f)=0$
for every $f \in F^* \backslash \spanv(\{f_1,\dots,f_p\})$. 
Then for every $i \in \{1, \dots,m\}$, $R(Z_i,\Phi)$. 
Consider the filter $\mathcal F$ on $I$ generated by $\{R[Z] ; Z \in fin(F^*)\}$. 
Then the mapping  $\Phi : F^* \to L(E,\IK_{\mathcal F})$
 associating to each $f \in F^*$ the $\IK$-linear mapping 
$\Phi(f):  E \to {\IK}_{\mathcal F}$  associating to each  
$x \in E$ the class of  $(f(x))_{f \in  I}$  in  $ {\IK}_{\mathcal F}$
 is linear, and for every $f \in F^*$, $\Phi(f) :E \to {\IK}_{\mathcal F}$ extends $f$. 
 Using $\D_{\IK}$, there exists (see Proposition~\ref{prop:D_red_puis}) a  $\IK$-linear retraction 
$r: {\IK}_{\mathcal F} \to \IK$ of 
$j_{\mathcal F}: \IK \to {\IK}_{\mathcal F}$. Let $T: F^* \to E^*$ be the mapping 
$f \mapsto r \circ {\Phi(f)}$. 
Then the mapping $T: F^* \to E^*$ is a linear extender for $F$ in $E$.
\Rp

\begin{remark}
Notice that the axiom $\LE_{\IK}$ is ``multiple'': given a family $(E_i)_{i \in I}$ of 
$\IK$-vector spaces and a family $(F_i)_{i \in I}$ such that for each $i \in I$, $F_i$ is a vector 
subspace of $E_i$, then there exists a family $(T_i)_{i \in I}$ such that for each $i \in I$,
$T_i: F_i^* \to E_i^*$ is a linear extender for $F_i$ in $E_i$: apply $\LE_{\IK}$ to the 
vector subspace $\oplus_{i \in I}F_i$ of $\oplus_{i \in I} E_i$.
\end{remark}

\begin{corollary} \label{coro:retrac}
Given a commutative field $\IK$, then $\D_{\IK}$ implies (and is equivalent to) the following 
statement: for every $\IK$-vector space $E$, 
for every vector subspace $F$ of $E$, denoting by $\can: F \to E$ the canonical mapping, then the 
double transposed mapping  $\can^{tt}: F^{**} \to E^{**}$ is one-to-one and has a $\IK$-linear 
retraction $r:E^{**} \to F^{**}$. 
\end{corollary}
\Pr The mapping $\can^{tt}: F^{**} \to E^{**}$ associates to every $\Phi \in F^{**}$ the 
mapping $\overline{\Phi}: E^{*} \to \IK$ such that for every $g \in E^*$,  
 $\overline{\Phi}(g)=\Phi(g_{\restriction F})$. 
Given some $\Phi \in \ker(\can^{tt})$, then,  for every 
$g \in E^*$,  $(\can^{tt}(\Phi))(g)=0$ {\em i.e.}  
$\Phi (g_{\restriction F})=0$; using $\D_{\IK}$, for every $f \in F^*$ there exists  
some $g \in E^*$ such that $f=g_{\restriction F}$, thus for every 
$f \in F^*$, $\Phi (f)=0$, so $\Phi=0$. It follows that $can^{tt}: F^{**} \to E^{**}$ is 
one-to-one. 
Using the equivalent form $\LE_{\IK}$ of $\D_{\IK}$,
let $T: F^* \to E^*$ be a $\IK$-linear extender {\em i.e.} a $\IK$-linear mapping such that for 
each 
$f \in F^*$, $T(f): E \to \IK$
extends $f$. Then the transposed mapping $T^t: E^{**} \to F^{**}$ is $\IK$-linear and
for every $\Phi \in F^{**}$, $T^t(\can^{tt}(\Phi))=\can^{tt}(\Phi)  \circ T = \Phi$. 
\Rp 

It follows that $\D_{\IK}$ implies  (and is equivalent to) the following statement: 
``For every vector space $E$, for every vector subspace $F$ of of $E$, the canonical linear mapping 
$F^{**} \rightarrow E^{**}$ is one-to-one and there exists a $\IK$-vector space $G$ such that  $G$ 
is a complement of $F^{**}$ in $E^{**}$. 

\begin{remark}
Corollary~\ref{coro:retrac} is equivalent to its ``multiple form'': given a family 
$(E_i)_{i \in I}$ of $\IK$-vector spaces and a family $(F_i)_{i \in I}$ such that for each $i \in 
I$, 
$F_i$ is a vector subspace of $E_i$, then for each $i \in I$, the canonical mapping
$\can_i: F_i^{**} \to E_i^{**}$ is one-to-one  and there exists a family $(r_i)_{i \in I}$ 
such that for each $i \in I$,
$r_i: E_i^{**} \to F_i^{**}$ is a $\IK$-linear retraction of $\can_i:F_i^{**} \to E_i^{**}$. 
\end{remark}

\section{Other statements equivalent to $\D_{\IK}$} \label{sec:other_equ}
\subsection{$\IK$-linearity and additive retractions}
\begin{proposition} \label{prop:add-to-Klinear2}
Let $\IK$ be a commutative  field and let $a$ be a non-null element of a $\IK$-vector space 
$E$. Let $j_a:\IK \hookrightarrow E$ be the mapping  $\lambda \mapsto \lambda.a$.
Given any additive mapping  $r: E \to \IK$ such that $r(a)=1$,
then $r$ is $\IK$-linear if and only if $r$ is a retraction of the  mapping   
$j_a:\IK \hookrightarrow E$.
\end{proposition}
\Pr  The direct implication is easy to prove. We show the converse statement. Assuming that 
$r: E \to \IK$ is an additive retraction of $j_a$, let us check that $r$ is $\IK$-linear. 
Since $r$ is a retraction of $j_a$, $\ker(r) \cap  \IK.a=\{0\}$. Thus 
 $\ker(r) \oplus \IK.a=E$ is the direct sum of groups with 
the unique decomposition $x=(x-r(x).a)+r(x).a$ for every $x \in E$. Therefore, 
$\ker(r)$ is a maximal subgroup $H$ of $E$ such that $H \cap \IK.a = \{0\}$. Also notice that, from 
$\ker(r) \cap  \IK.a=\{0\}$, it follows that $\IK.\ker(r) \cap \IK.a=\{0\}$, thus 
$\IK.\ker(r)=\ker(r)$.  We now check that  the additive mapping $r$ is $\IK$-linear: given $z \in 
E$, then $z=x \oplus t.a$ where $x \in \ker(r)$ and 
$t \in \IK$; thus for every $\lambda \in \IK$, $r(\lambda.z)=r(\lambda.x \oplus \lambda t.a)
=r(\lambda.x) + r(\lambda t.a)=r(\lambda t.a)=\lambda t=\lambda.r(z)$. 
\Rp

\subsection{Additivity  statements equivalent to  $ \D_{\IK}$}
Given a commutative field $(\IK,+,\times,0,1)$, we  consider the following statements:
\par $\A_{\IK}$: {\em ``For every  $\IK$-vector space $E$ and every subgroup  
$F$ of $(E,+)$, for every additive  mapping $f: F \to (\IK,+)$, there exists an  additive mapping 
$\tilde f:E \to \IK$  extending $f$.}
\par $\A'_{\IK}$: {\em ``For every $\IK$-vector space $E$ and every vector subspace 
$F$ of $E$, for every $\IK$-linear  mapping $f: F \to \IK$, there exists an  additive mapping 
$\tilde f:E \to \IK$ extending $f$.}
\par $\A''_{\IK}$: {\em ``For every  $\IK$-vector space $E$ and every 
$a \in E \backslash \{0\}$, there exists an additive mapping $f: E \to \IK$ 
such that for every $\lambda \in \IK$, $f(\lambda.a)=\lambda$.''}

\begin{proposition}
For every commutative field $\IK$, 
$  \A_{\IK}  \Leftrightarrow \A'_{\IK} \Leftrightarrow \A''_{\IK} \Leftrightarrow \D_{\IK}$.
\end{proposition}
\Pr The implications $ \A_{\IK} \Rightarrow \A'_{\IK}$ and $ \A'_{\IK} \Rightarrow \A''_{\IK}$
are trivial. We prove $\A''_{\IK} \Rightarrow \D_{\IK}$. 
In view of Proposition~\ref{prop:D_red_puis}, we prove that for every filter $\mathcal F$ on a set 
$I$, 
the canonical mapping 
$j_{\mathcal F} : \IK \hookrightarrow  \IK_{\mathcal F}:={\IK}^I/{\mathcal F}$ has a $\IK$-linear 
retraction.  Let $f: \IK \to \IK$ be the identity mapping. 
Using $\A''_{\IK}$, let $\tilde f: \IK_{\mathcal F} \to \IK$ be an additive mapping extending 
$f$. 
Using Proposition~\ref{prop:add-to-Klinear2},  $\tilde f$ is $\IK$-linear. It follows that $\tilde f$
is a $\IK$-linear retraction of $j_{\mathcal F} : \IK \hookrightarrow \IK_{\mathcal F}$. \\
$ \D_{\IK} \Rightarrow \A_{\IK}$. Let $E$ be a $\IK$-vector space, let $F$ be a subgroup of the 
additive group $(E,+)$, and let $f: F \to (\IK,+)$ be an additive mapping. Using a concurrent 
relation (the proof is similar to the proof of  Lemma~1 in \cite{Mo09}), let $\mathcal F$ be a 
filter on the set $I=\IK^E$ and let 
$\iota: E \to \IK_{\mathcal F}$ be an additive mapping extending $f$.  Using $\D_{\IK}$,
let $r: \IK_{\mathcal F} \to \IK$ be an additive retraction of 
$j_{\mathcal F} : \IK \to \IK_{\mathcal F}$. Then $\tilde f:= r \circ \iota: E \to \IK$ is additive 
and extends $f$.
\Rp

\subsection{A consequence of $\D_{\IK}$}
\begin{proposition}
Given a commutative field $\IK$, the following statements are equivalent:
\begin{enumerate}[i)]
\item \label{it:wD1}  $w\D_{\IK}$: {\em ``For every  $\IK$-vector space $E$ and every 
$a \in E \backslash \{0\}$, there exists an additive mapping $f: E \to \IK$ 
such that  $f(a)=1$.''} 
\item  \label{it:wD2} {\em ``For every  non null $\IK$-vector space $E$, there exists a 
non null additive mapping $f: E \to \IK$.''}
\item  \label{it:wD3} {\em ``For every  filter $\mathcal F$ on a set $I$, there exists a non null 
additive mapping $f: {\IK}_{\mathcal F} \to \IK$.''}
\end{enumerate}
\end{proposition}
\Pr  \eqref{it:wD1} $\Rightarrow$ \eqref{it:wD2}  and  \eqref{it:wD2} $\Rightarrow$ \eqref{it:wD3}  
are easy. We prove  \eqref{it:wD3} $\Rightarrow$ \eqref{it:wD1} Given a non null element $a$ of a 
$\IK$-vector space $E$, let $\mathcal F$ be a filter on a set $I$ and $g: E \to \IK_{\mathcal F}$ 
be 
a $\IK$-linear mapping such that $g(a)=1$. Using \eqref{it:wD3}, let 
$r: \IK_{\mathcal F} \to \IK$  be a non null additive mapping. Let $\alpha \in \IK_{\mathcal F}$
such that $r(\alpha) \neq 0$. Let $m_{\alpha} : \IK_{\mathcal F} \to  \IK_{\mathcal F}$ be the
additive mapping $x \mapsto \alpha x$. Then $g_1:=r \circ m_{\alpha} \circ g: E \to \IK$ is 
additive; let $f:=\frac{1}{r(\alpha)}.g_1 : E \to \IK$; then $f$ is additive and $f(a)=1$.  
\Rp

\begin{remark} Given a commutative field $\IK$ with prime field $k$, then $\D_k$ implies 
$w\D_{\IK}$.  
\end{remark}
\Pr Given a $\IK$-vector space $E$, a mapping $f: E \to \IK$ is additive if and only if
$f$ is $k$-linear.
\Rp

\begin{question}
Given a commutative field $\IK$ with null characteristic, does $\D_{\IQ}$ imply $\D_{\IK}$?
\end{question}

\begin{question}
Given a prime number $p$ and a commutative field $\IK$ with characteristic $p$, and denoting by 
$\IF_p$ the finite field with $p$ elements, does $\D_{\IF_p}$ imply $\D_{\IK}$? 
Does $\BPI$ (which implies $\D_{\IF_p}$) imply $\D_{\IK}$?
\end{question}

\section{Ingleton's statement for ultrametric  valued fields} \label{sec:ingleton}
\subsection{Semi-norms on vector spaces over a valued field}
\subsubsection{Pseudo metric spaces}
Given a  set $X$, a {\em pseudo-metric} on $X$ is a mapping $d: X \times X \to \IR_+$ 
such that for every $x,y,z \in X$, $d(x,y)=d(y,x)$ and $d(x,z) \le d(x,y)+d(y,z)$. If $d $ satisfies 
the extra property $(d(x,y)=0 \Rightarrow x=y)$, then $d$ is a {\em metric} on $X$.
A pseudo-metric $d$ on $X$ is said to be {\em ultrametric} if for every $x,y,z \in X$,
$d(x,z) \le \max(d(x,y), d(y,z))$. 

Given a pseudo-metric space $(X,d)$, for every $a \in X$ and every $r \in \IR_+$,  
we denote by $B_s(a,r)$ the ``strict'' ball
$\{x \in X : d(x,a)<r\}$ and we denote by $B(a,r)$ the ``large'' ball
$\{x \in X : d(x,a) \le r\}$. Notice that large balls of a pseudo-metric space are nonempty.
A pseudo-metric space $(X,d)$ is {\em spherically complete} if every chain ({\em i.e.} 
set which is linearly ordered for the inclusion) of  large balls of $X$ has a nonempty intersection.

\begin{example}
Given a nonempty set $X$, the discrete metric $d_{disc}$ on $X$, associating to each 
$(x,y) \in X \times X$ the real number 1 if $x \neq y$ and 0 else is ultrametric and the 
associated metric space $(X,d_{disc})$ is spherically complete since large balls for this metric  
are singletons and the whole space $X$.  
\end{example}

\subsubsection{Group semi-norms}
Given a  commutative group  $(G,+,0)$, a {\em group semi-norm}  on $G$ is a  mapping 
$N: G \to \IR_+$ which is  sub-additive  (for every  $x,y \in G$, $N(x+y) \le N(x)+N(y)$) 
and symmetric (for every $x \in G$, $N(-x)=N(x)$) and such that $N(0)=0$. 
If for every $x \in G$, $(N(x)=0 \Rightarrow x=0)$, then $N$ is  a {\em norm}. 
Given a group semi-norm $N$ on an abelian group $(G,+,0)$, the mapping $d: G \times G \to \IR_+$ 
associating to each $(x,y) \in G \times G$ the real number $N(x-y)$ is a pseudo-metric on  $G$. 
Moreover, if $N$ is a norm, then  $d$ is a metric on $G$. The semi-norm $N$ is said to be 
ultrametric if the pseudo-metric $d_N$ is ultrametric, which is equivalent to say that for every 
$x,y \in G$, $N(x+y) \le \max(N(x),N(y))$.
 The topology on $G$ associated to the pseudo-metric  $d_N$ is a group topology on $G$, which is 
 Hausdorff if and only if $N$ is a norm. 

\subsubsection{Absolute values on a commutative field}
Given a  commutative unitary ring  $(R,+,\times,0,1)$, a {\em ring semi-norm} 
(see \cite[p.~137]{warner}) on $R$ is a 
group semi-norm on $(R,+,0)$ which is sub-multiplicative: for every $x,y \in R$, 
$N(x \times y) \le N(x)N(y)$; if a ring semi-norm $N:R \to \IR_+$ is non null, then $N(1) \ge 1$. 
  If $R$ is a commutative field and if $N$ is a ring semi-norm on $R$ which is multiplicative 
  (for every $x,y \in R$, $N(x \times y) = N(x)N(y)$) and non null, then $N(1)=1$,  $N$ is a norm 
and $N$ is called an {\em absolute value} on the commutative field $R$, and 
$(R,N)$ is a {\em valued field}.

\begin{example}
Given a  commutative field  $\IK$,  the mapping $|.|_{disc}: \IK \to \IR_+$ associating to 
each $x \in \IK$ the real number $1$ if $x \neq 0$ and $0$ else is  an absolute value on $\IK$, 
which is called the {\em trivial} absolute value on $\IK$. The metric associated to this absolute 
value is the discrete metric $d_{disc}$ on $\IK$, thus, the discrete field $(\IK,|.|_{disc})$ is 
spherically complete.
\end{example}

\subsubsection{Vector semi-norms  on vector spaces over a valued field}
Given a  valued field $(\IK,|.|)$  and a $\IK$-vector space  $E$, a {\em vector semi-norm} on 
$E$  (see \cite[p.~210]{warner}) is a group semi-norm 
$p: E \to \IR_+$ on the abelian group $(E,+)$  such that for every 
$\lambda \in \IK$ and every $x \in E$, $p(\lambda.x)=|\lambda| p(x)$; if $(\IK,|.|)$ is 
ultrametric, 
then the semi-norm $p: E \to \IR_+$ is said to be {\em ultrametric} if and only if for every $x,y 
\in E$, $p(x+y) \le \max(p(x),p(y))$. 
Given a  $\IK$-algebra $(A,+,\times,0,1,\lambda .)$, an {\em algebra semi-norm} on $A$
is a ring semi-norm on the ring $A$ which is also a vector semi-norm on the $\IK$-vector space $A$. 

\subsection{Ingleton's statement for spherically complete ultrametric valued fields} 
\label{subsec:ingleton}
Given a valued field $(\IK,|.|)$, a set $E$ and a mapping $f:E \to \IK$, we denote by 
$|f| : E \to \IR_+$ the mapping $x \mapsto |f(x)|$. We endow $\IR^E$ with the product order, thus 
given mappings $p: E \to \IR$ and $q: E \to \IR$, $p \le q$ means that for every $x \in E$, $p(x) 
\le q(x)$. In particular, $|f| \le p$ means that for every $x \in E$, $|f(x)| \le p(x)$.

\begin{lemma*}[Ingleton, \cite{Ing52}]
Let $(\IK,|.|)$ be a  ultrametric spherically complete valued  field, let $E$ be a $\IK$-vector 
space and 
let  $p: E \to \IR$ be an ultrametric vector semi-norm. Let $F$ be a vector subspace of $E$ and let 
$f: F \to \IK$ be a linear mapping such that $|f| \le p$. 
Then, for every $a \in E \backslash F$, there exists 
a linear mapping $\tilde f: F \oplus \IK a \to \IK$ extending $f$ such that 
$\tilde f \le p_{\restriction F \oplus \IK a}$.
\end{lemma*}
\Pr Ingleton's proof of this Lemma  holds in $\ZF$. For sake of completeness, we give a proof of 
this Lemma. We search for some $\alpha \in \IK$ such that:
\begin{equation} \label{it:ing-alpha}
\forall z \in F \; \forall \lambda \in \IK \; |f(z)+\lambda. \alpha| \le p(z+\lambda.a)
\end{equation}
Given some $\alpha \in \IK$, and denoting by $\IK^*$ the set $\IK \backslash \{0\}$, 
Condition~\eqref{it:ing-alpha} is equivalent to:
$$ \forall z \in F \;  \forall \lambda \in \IK^* \;
 |\lambda. (f(\frac{z}{\lambda})+ \alpha)| \le p(\lambda(\frac{z}{\lambda}+a)) $$
$$ \forall z \in F \;  \forall \lambda \in \IK^*  \;
 |\lambda|. |f(\frac{z}{\lambda})+ \alpha)| \le |\lambda|. p(\frac{z}{\lambda}+a) $$
$$ \forall z \in F \;  \forall \lambda \in \IK^* \;
 |(f(\frac{z}{\lambda})+ \alpha)| \le p(\frac{z}{\lambda}+a) $$
$$ \forall z' \in F \;  |f(z')+ \alpha)| \le p(z'+a) $$
$$ \forall z \in F \;  |f(z)- \alpha)| \le p(z-a) $$
Let $\mathcal B_a$ the set of large balls $B(f(z),p(z-a))$ for $z \in F$.
The element  $\alpha \in \IK$  satisfies \eqref{it:ing-alpha} if and only if  
 $\alpha \in \cap  \mathcal B_a$.
 Given two  elements $z_1,z_2 \in F$,  then 
$|f(z_1) - f(z_2)| = |f(z_1 -z_2)| \le p(z_1-z_2) \le \max(p(z_1-a), p(z_2-a))$ since $p$ is 
ultrametric. If $p(z_1-a) \le p(z_2-a)$,
 it follows that $f(z_1) \in B(f(z_2),p(z_2-a))$
thus, since $(\IK,|.|)$ is ultrametric,  $B(f(z_1),p(z_2-a)) \subseteq B(f(z_2),p(z_2-a))$, 
thus $B(f(z_1),p(z_1-a)) \subseteq B(f(z_2),p(z_2-a))$, so  the set of large balls 
$\mathcal B_a$ is a chain. Since $(\IK,|.|)$ is spherically complete, it follows that 
$\cap \mathcal B_a$ is nonempty.
\Rp  

Given a spherically complete ultrametric valued field $(\IK,|.|)$, we now consider Ingleton's 
statement:\\
$\I_{(\IK,|.|)}$: {\em ``For every 
$\IK$-vector space $E$  endowed with an ultrametric vector semi-norm $p:E \to \IR_+$, for every 
vector subspace   $F$ of $E$ and every linear mapping $f: F \to \IK$ such that 
$|f| \le p_{\restriction F}$, there exists a  linear mapping $\tilde f: E \to \IK$ extending $f$  
such that
 $|\tilde f|  \le p$.''}

\begin{corollary*}[Ingleton,  \cite{Ing52}]
For every  spherically complete ultrametric valued field $(\IK,|.|)$, 
$\AC \Rightarrow \I_{(\IK,|.|)}$. 
\end{corollary*}
\Pr Use Zorn's lemma and the previous Lemma. 
\Rp

\begin{remark} \label{rem:ing_basis}
Notice that in $\ZF$, given a spherically complete ultrametric valued field $(\IK,|.|)$ and  
a $\IK$-vector space $E$ which has a well-orderable basis (for example a finitely generated 
$\IK$-vector space), 
then every ultrametric vector semi-norm $p:E \to \IR_+$ satisfies the Ingleton statement.
\end{remark}

\subsection{$\MC$ implies Ingleton's statement for null characteristic fields} 
\begin{proposition}
Given a spherically complete ultrametric valued field $(\IK,|.|)$ with null characteristic such that the restricted absolute value $|.|_{\restriction \IQ}$ is trivial,  
$\MC$ implies $\I_{(\IK,|.|)}$. It follows that in $\ZFA$, the Ingleton statement 
restricted to null characteristic ultrametric valued fields $(\IK,|.|)$ such that the restricted absolute value $|.|_{\restriction \IQ}$ is trivial does not imply $\AC$.
\end{proposition}
\Pr   Assume that $E$ is a $\IK$-vector space endowed with an ultrametric semi-norm 
$p:E \to \IR_+$,  and assume that $F$ is a vector subspace of $E$ and that $f: F \to \IK$ is 
a linear mapping such that $|f| \le p$. 
Using $\MC$, there exists an ordinal
$\alpha$ and some partition $(F_i)_{i \in  \alpha}$ in finite sets of $E \backslash F$. 
This implies that there is an ordinal $\beta$ and a strictly increasing family $(V_i)_{i \in  
\beta}$ of vector subspaces of $E$ such that $V_0=F$,  for every $i \in \beta$ such that $i+1 \in 
\beta$,
$V_{i+1}/V_i$ is finite-dimensional, and for every non null limit ordinal $i \in \beta$,
$V_i=\cup_{j <i}V_j$. 
Let $Z$ be the set of 3-uples $(V,W,l)$ such that $V,W$ are  subspaces of $E$ satisfying 
$V \subseteq W$, 
$W/V$ is finite-dimensional and  $l:V \to \IK$ is a  linear mapping satisfying $|l| \le p$.
For each $(V,W,l) \in Z$, using the previous Lemma  (which holds in $\ZFA$), the set 
$A_{(V,W,l)}$ of linear mappings $\tilde l: W \to \IK$ extending $l$ and satisfying $|\tilde l|\le 
p$ 
is non-empty; using $\MC$, there is a mapping 
associating to every $i=(V,W,l) \in Z$ a nonempty finite subset  
$B_i$ of $A_{(V,W,l)}$. 
Then, for every $i \in Z$, define $\Phi(i) := \frac{1}{\#B_i} \sum_{u \in B_i} u$ where $\#B_i$ is the cardinal of the nonempty finite set $B_i$ 
(here we use the fact that the characteristic of $\IK$ is null); notice that $\Phi(i)$ is linear and that 
that for every $x \in W$, $|\Phi(i)(x)|=|\sum_{u \in B_i} u(x)| \le \max_{u \in B_i} |u(x)| \le p(x)$. 
 Using the function $\Phi$, we define by transfinite recursion a family 
$(f_i)_{i \in  \beta}$ such that for each $i \in  \beta$, $f_i : V_i \to \IK$ is linear, 
$f_i$ extends $f$, $|f_i|\le p$, 
 and for every
$i<j \in \beta$, $f_j$ extends $f_i$. Let $\tilde f:= \cup_{i \in \beta} f_i$. Then 
$\tilde f: E \to \IK$ is linear, $\tilde f$ extends  $f$ and $|\tilde f| \le p$. 
\Rp

\begin{question}
Given a spherically complete ultrametric valued field $(\IK,|.|)$, 
does $\B_{\IK}$ imply $\I_{(\IK,|.|)}$?  Does $\BE_{\IK}$ imply 
$\I_{(\IK,|.|)}$?
\end{question}

\begin{question}[van Rooij, see \cite{Rooij}]
Does the full Ingleton statement
 ({\em i.e.} $\I_{(\IK,|.|)}$ for every spherically complete ultrametric valued 
 field $(\IK,|.|)$) imply $\AC$? 
\end{question}

\begin{remark}
If a commutative field $\IK$ is endowed with the trivial absolute value 
$|.|_{disc}$,  then $\I_{(\IK,|.|_{disc})}$ is the statement $\D_{\IK}$, thus 
the ``full Ingleton statement'' implies $\D_{\IK}$ for every commutative field $\IK$. 
\end{remark}
\Pr $\I_{(\IK,|.|_{disc})} \Rightarrow \D_{\IK}$. Given a $\IK$-vector space $E$ and a non null 
vector $a \in E$, consider the trivial ultrametric vector semi-norm $p: E \to \IR_+$ defined by  
$p(0_E)=0$ and 
for every $x \in E \backslash \{0\}$, $p(x)=1$. Let $f: \IK.a \to \IK$ be the linear mapping such 
that $f(a)=1$; then $|f(0_E)|_{disc}=0=p(0_E)$ and, for every $\lambda \in \IK \backslash \{0\}$, 
$|f(\lambda.a)|_{disc} =|\lambda|_{disc}=1=p(\lambda.a)$. Using $\I_{(\IK,|.|_{disc})}$, there 
exists a linear mapping $\tilde f: E \to \IK$ extending $f$ such that for every $x \in E$,
$|\tilde f(x)|_{disc} \le p(x)$. Thus $\tilde f$ is a linear form on $E$ such that
$\tilde f(a)=1$.\\
$\D_{\IK} \Rightarrow \I_{(\IK,|.|_{disc})}$. Let $E$ be a $\IK$-vector space and let $p: E \to 
\IR_+$ be an ultrametric vector semi-norm with respect to the valued field $(\IK,|.|_{disc})$. 
Assume that 
$F$ is a vector subspace of $E$ and that $f: F \to \IK$ is a linear form such that for every 
$x \in F$, $|f(x)|_{disc} \le p(x)$. If $f$ is null, then the null mapping $\tilde f: E \to \IK$ 
extends $f$ and satisfies $|\tilde f(x)|_{disc} \le p(x)$ for every $x \in E$. If $f$ is not null, 
let $a \in F$ such that $f(a)=1$. Denoting by $V$ the vector subspace $\{x \in E : 
p(x)<1\}$, then $a \notin V+\ker(f)$ because  if $a=a_1+a_2$ with $a_1 \in V$ and $a_2 \in 
\ker(f)$, then $a_1=a-a_2 \in F$ thus 
$1=f(a)=f(a_1) \le p(a_1)$ so  $a_1 \notin V$, which is contradictory!
Let $\can: E \to E/(V+\ker(f))$ be the quotient mapping.  Since $a \notin V+\ker(f)$, it 
follows that  $\can(a) \neq 0$.
Using $\D_{\IK}$,  let $g: E/(V+\ker(f)) \to \IK$ be a linear mapping such that $g(\can(a))=1$. 
Then $\tilde f:=g \circ \can: E \to \IK$ is a linear mapping which is null on $V+\ker(f)$  and 
such 
that $\tilde f(a)=1$. Since $F=\ker(f) \oplus \IK.a$, it follows that $\tilde f$ extends $f$.
We now show that $|\tilde f | \le p$. Given some  $x \in E$, if $x \in V$ then 
$\tilde f(x)=0$ 
so $|\tilde f(x)|_{disc} \le p(x)$; else  $p(x) \ge 1$, thus $|\tilde f(x)|_{disc} \le 1 \le 
p(x)$.
\Rp

\subsection{Ingleton's statement for ultrametric fields with compact large balls}
Given a valued field $(\IK,|.|)$ and a filter $\mathcal F$ on a set $I$, we denote by 
$|.|_{\mathcal F}: \IK_{\mathcal F} \to \IR_{\mathcal F}$ the quotient mapping
associating to each $x \in \IK_{\mathcal F}$ which is the equivalence class of some 
$(x_i)_{i \in I} \in \IK^I$, the class of $(|x_i|)_{i \in I}$ in $\IR_{\mathcal F}$.
We denote by $(\IK_{\mathcal F})_b$ the following unitary subalgebra of  ``bounded elements'' of  
the unitary  $\IK$-algebra $\IK_{\mathcal F}$:
$\{x \in \IK_{\mathcal F} : \exists t \in \IR_+ : |x|_{\mathcal F} \le_{\mathcal F} t\}$. 
We also denote by 
$N_{|.|,\mathcal F}: (\IK_{\mathcal F})_b \to \IR_+$ the mapping associating to each
$x \in  (\IK_{\mathcal F})_b$ the real number 
$\inf \{t \in \IR_+ : |x|_{\mathcal F} \le_{\mathcal F} t\}$. The mapping 
$N_{|.|,\mathcal F}$ is a unitary algebra  semi-norm on $(\IK_{\mathcal F})_b$; 
moreover, if the valued field $(\IK,|.|)$ is ultrametric, the vector  semi-norm $N_{|.|,\mathcal F}$ 
is also ultrametric.

\begin{lemma} \label{lem:ext-valeurs-ns}
Let  $(\IK,|.|)$ be a spherically complete ultrametric valued field, let $E$ be a  $\IK$-vector 
space and let $p: E \to \IK$ be an ultrametric vector semi-norm. Assume that $F$ is a vector 
subspace 
of $E$ and that 
$f: F \to \IK$ is a linear form such that $|f| \le p$. Let  $I:=\IK^E$. There exists   a filter 
$\mathcal F$  on  $I$ and a $\IK$-linear mapping $\iota: E \to (\IK_{\mathcal F})_b$ definable 
from $E$, $p$ and $f$ such that
$\iota$ extends $f$ and such that $N_{|.|,\mathcal F} \circ \iota \le p$.
\end{lemma}
\Pr Let $R \subseteq (fin(E) \times I)$ be the following binary relation: given $Z \in fin(E)$ and 
given some mapping $u : E \to \IK$, then $R(Z,u)$  iff $u$ extends $f$, $|u| \le p$ and 
$u_{\restriction Z}$ is linear {\em i.e.}  for every $x,y \in Z$ and $\lambda \in \IK$, 
$(x+y \in Z \Rightarrow u(x+y)=u(x)+u(y))$ and 
$(\lambda x \in Z \Rightarrow u(\lambda x)=\lambda u(x))$. Using Ingleton's Lemma in 
Section~\ref{subsec:ingleton}, the binary relation $R$ is concurrent, thus it generates a filter
$\mathcal F$ on $I$.  
Let  $\iota : E \to \IK_{\mathcal F}$ be the mapping associating to each $x \in E$, 
the equivalence class of $(i(x))_{i \in I}$ in $ \IK_{\mathcal F}$. 
Then, for every $x \in F$,   $\iota(x)=f(x)$ and for every $x \in E$, 
$|\iota(x)|_{\mathcal F} \le p(x)$ whence $N_{|.|,\mathcal F}(\iota(x)) \le p(x)$. 
Moreover, $\iota$ is $\IK$-linear: given $x, y \in E$ and 
$\lambda \in \IK$,  let $Z:=\{x,y, \lambda y, x+\lambda y\}$; by definition of $\iota$,  
the set 
$J:=\{i \in I: R(Z,i)\}$ belongs to 
$\mathcal F$, and 
$J \subseteq \{i \in I: i(x+\lambda y)=i(x)+ \lambda i(y)\}$; thus 
$\iota(x+ \lambda y)=\iota(x)+ \lambda \iota(y)$ so $\iota$ is $\IK$-linear.
\Rp

\begin{lemma} \label{lem:equiv-ing}
Given a spherically complete ultrametric valued field $(\IK,|.|)$, the following statements are 
equivalent:
\begin{enumerate}[(i)]
\item \label{it:ax-ing1} $\I_{(\IK,|.|)}$
\item  \label{it:ax-ing2} For every filter $\mathcal F$ on a set $I$, 
the $\IK$-linear mapping $j_{\mathcal F} : \IK \to (\IK_{\mathcal F})_b$ has an additive
retraction $r: (\IK_{\mathcal F})_b \to \IK$  such that for every 
$x \in  (\IK_{\mathcal F})_b$, $|r(x)| \le N_{|.|,\mathcal F}(x)$.
\end{enumerate}
\end{lemma}
\Pr \eqref{it:ax-ing1} $\Rightarrow$  \eqref{it:ax-ing2}  Given a filter $\mathcal F$ on a set 
$I$, consider the $\IK$-linear mapping 
$f: \IK.1_{\mathcal F} \to \IK$ associating $1$ to $1_{\mathcal F}$. Then 
$|f| \le p_{\mathcal F}$. Since the vector semi-norm
$p_{\mathcal F}: (\IK_{\mathcal F})_b \to \IR_+$ is  ultrametric,
\eqref{it:ax-ing2} implies a $\IK$-linear mapping $r:  (\IK_{\mathcal F})_b \to \IK$
extending $f$ such that $r \le N_{|.|,\mathcal F}$. Since $r$ is $\IK$-linear and fixes 
$1_{\mathcal F}$, $r$ fixes every element of $\IK$ in $\IK_{\mathcal F}$ thus $r$ is a retraction 
of $j_{\mathcal F} : \IK \to (\IK_{\mathcal F})_b$. \\
 \eqref{it:ax-ing2} $\Rightarrow$  \eqref{it:ax-ing1} Given an ultrametric vector semi-norm 
$p$ on a $\IK$-vector space $E$, and a linear mapping $f$ defined on a vector subspace $F$ of 
$E$ such that $|f| \le p$, using 
Lemma~\ref{lem:ext-valeurs-ns}, consider a linear mapping $\iota: E \to  (\IK_{\mathcal F})_b$ 
extending $f$ such that $N_{|.|,\mathcal F} \circ \iota \le p$. Using \eqref{it:ax-ing2}, let
$r: \IK_{\mathcal F} \to \IK$ be an additive  retraction such that 
$|r| \le N_{|.|,\mathcal F}$; using Proposition~\ref{prop:add-to-Klinear2}, $r$ is $\IK$-linear. 
Then the $\IK$-linear mapping  $\tilde f:=r \circ \iota: E \to \IK$ 
extends $f$ and $|\tilde f| \le p$. 
\Rp

\begin{remark}[Hahn-Banach]
Consider the Hahn-Banach statement $\HB$: ``Given a vector space $E$ over $\IR$,
given a subadditive mapping $p: E \to \IR$ such that for every $\lambda \in \IR_+$ and
every $x \in E$, $p(\lambda.x)=\lambda p(x)$, and given a linear form $f$ defined on a vector 
subspace $F$ of $E$ satisfying $|f| \le p$, there exists a linear mapping $\tilde f: E \to \IR$
extending $f$ such that $|\tilde f| \le p$''. It is known (see \cite {Ho-Ru}) 
that $\BPI \Rightarrow \HB$ and that 
$\HB \not \Rightarrow \BPI$.  It is also known (see \cite{Lu})  that $\HB$ is equivalent to 
the following statement: ``For every filter $\mathcal F$ on a set $I$, there exists a $\IR$-linear 
mapping
$r: (\IR_{\mathcal F})_b \to \IR$ such that $r(1)=1$ and $r$ is positive.
\end{remark}

\begin{question} Is there an ultrametric spherically complete valued field $(\IK,|.|)$ such that 
$\HB$ is equivalent to $\I_{(\IK,|.|)}$? Given two distinct prime numbers $p$ and $q$, 
are the statements $\I_{\IQ_p}$ and $\I_{\IQ_q}$ equivalent? Are they equivalent to $\HB$? Here we 
denote by $\IQ_p$   (see \cite[p.~186]{warner}) the valued field which is the Cauchy-completion of 
$\IQ$ endowed with the $p$-adic absolute value.
\end{question}

\begin{remark}
Every  ultrametric valued field  in which every large ball is compact is spherically complete.
\end{remark}

\begin{corollary*}[van Rooij, \cite{Rooij}]
For every   ultrametric valued field  
$(\IK,|.|)$ such that every large ball of $\IK$ is compact, then $\BPI$ implies $\I_{(\IK,|.|)}$.
\end{corollary*}
\Pr For sake of completeness, we give the proof  sketched by van Rooij. 
Using  Lemma~\ref{lem:equiv-ing}, it is sufficient to show that given a filter $\mathcal F$ on a 
set $I$, there is a $\IK$-linear mapping $r: (\IK_{\mathcal F})_b \to \IK$ such that 
$|r| \le N_{|.|,,\mathcal F}$.   Using 
$\BPI$, let $\mathcal U$ be an ultrafilter on $I$ such that  $\mathcal F \subseteq \mathcal U$. 
For every $x  \in (\IK_{\mathcal F})_b$, 
the large ball $B(0, N_{|.|,\mathcal F}(x)+1)$ of $\IK$ is compact and Hausdorff, whence for every 
$(u_i)_{i \in I} \in \IK^I$ and 
 $(v_i)_{i \in I} \in \IK^I$ such that $x$ is the class of $(u_i)_{i \in I}$ 
and the class of $(v_i)_{i \in I}$  in $\IK_{\mathcal F}$, then $(u_i)_{i \in I}$ and 
$(v_i)_{i \in I}$ both  converge through the ultrafilter $\mathcal U$ to the same element 
of the ball   $B(0, N_{|.|,\mathcal F}(x)+1)$:
we denote by $r(x)$ this element of $\IK$. Since the class of $(u_i)_{i \in I}$ in 
$\IK_{\mathcal F}$ is $x$, for every real number $\varepsilon>0$, $|r(x)| \le 
N_{|.|,,\mathcal F}(x)+\varepsilon$,
thus $|r(x)| \le N_{|.|,,\mathcal F}(x)$. 
We have defined a mapping $r: (\IK_{\mathcal F})_b \to \IK$.
 Then $r: (\IK_{\mathcal F})_b \to \IK$  is additive and fixes every element of $\IK$, thus,
 using Proposition~\ref{prop:add-to-Klinear2}, $r$ 
 is $\IK$-linear. And $|r| \le  N_{|.|,,\mathcal F}$ by construction.
\Rp

\begin{remark}
In particular, given a  finite field $\IK$ endowed with the trivial absolute value, then  $\IK$ is 
compact thus 
$\BPI$ implies $\D_{\IK}$: this is 
Howard and Tachtsis's result (see \cite[Theorem~3.14]{How-Tach13}). 
\end{remark}

\begin{question}[van Rooij, see \cite{Rooij}]
Does $\BPI$ imply the full Ingleton statement?
\end{question}

\section{Isometric linear extenders} \label{sec:cont-le}
\subsection{Bounded dual of a semi-normed vector space over a valued field}
Given a  valued field $(\IK,|.|)$ and two semi-normed $\IK$-vector spaces $(E,p)$
and $(F,q)$, a linear mapping $T: E \to F$ is said to be {\em bounded}  with respect to the 
 semi-norms $p$ and $q$ if and only if there exists a real number 
$M \in \IR_+$ satisfying $q(T(x)) \le M p(x)$ for every $x \in E$. 

\begin{proposition*}[{\cite[Proposition~3.1. p.~13]{Schneider}}]
Let $(\IK,|.|)$ be a valued field and let $(E,p)$ and $(F,q)$ be two semi-normed $\IK$-vector 
spaces. Let  $T: E \to F$ be a $\IK$-linear mapping. 
\begin{enumerate}
 \item \label{it:bound2cont} If $T$ is bounded with respect to the 
 semi-norms $p$ and $q$, then $T$ is continuous with respect to the 
topologies associated to the semi-norms $p$ and $q$.
\item \label{it:cont2bound}If $T$ is continuous with respect to the 
topologies associated to the semi-norms $p$ and $q$ and if  the absolute value $|.|$ 
on $\IK$ is not trivial, then $T$ is  bounded with respect to the semi-norms $p$ and $q$.
\end{enumerate}
\end{proposition*}
\Pr \eqref{it:bound2cont} If $T$ is bounded, then $T$ is continuous at the point $0_E$ of $E$. 
Since translations of $E$ are continuous with respect to $p$, it follows that $T$ is continuous at 
every point of $E$.\\
\eqref{it:cont2bound} We assume that the absolute value $|.|$ on $\IK$ is not trivial. Let $G$ be 
the  subgroup $\{|x| : x \in \IK \backslash \{0\}\}$ of $(\IR^*_+,\times)$. Since $|.|$ is not 
trivial, there exists $a \in \IK^*$ such that $|a| \neq 1$. Using $\frac{1}{a}$ instead of $a$, we 
may assume that $0 < |a| <1$. It follows that the sequence $(|a^n|)_{n \in \IN}$ of $\IR_+^*$ 
converges to $0$. Since $T$ is continuous at point $0_E$, let $\eta \in \IR^*_+$ such 
that for every $x \in E$, $(p(x)<\eta \Rightarrow q(T(x)) <1)$. Let $n_0 \in \IN$ such that 
$|a^{n_0}|<\eta$ and let $M:=\frac{1}{|a^{n_0+1}|}$.  Then for every $x \in E$, let us check that 
$q(T(x)) \le Mp(x)$. If $p(x)=0$, then for every $\lambda \in \IK^*$, $p(\lambda.x)=0$ hence 
$q(T(\lambda x))<1$, thus  for every $\lambda \in \IK^*$, $q(T(x))<\frac{1}{|\lambda|}$ so 
$q(T(x))=0$. If $p(x)>0$, let 
$n \in \IN$ such that $|a^{n+1}| \le p(x) < |a^{n}|$; then $p(\frac{x}{a^n})<1$ thus
$p(\frac{a^{n_0}x}{a^n})<|a^{n_0}|$ so $q(T(\frac{a^{n_0}x}{a^n}))<1$ {\em i.e.}
$q(T(x))<|a^{n-n_0}|=\frac{|a^{n+1}|}{|a^{n_0+1}|} \le \frac{p(x)}{|a^{n_0+1}|} =Mp(x)$.
\Rp

\begin{remark}[{\cite[Example~3 p.77-78]{Na-Be-valu}}]
Given a  commutative valued field $\IK$ endowed  with the trivial absolute value $|.|_{disc}$, and 
two semi-normed $\IK$-vector spaces $(E,p)$
and $(F,q)$,  a continuous linear mapping $T: E \to F$ is   not necessarily bounded with respect to 
the 
 semi-norms $p$ and $q$. For sake of completeness, we sketch the argument. 
Let $E$ be 
the ring $\IK[X]$ of polynomials with coefficients in $\IK$, let $p$ be the trivial norm on 
$\IK[X]$ and let $q: \IK[X] \to \IR_+$ be the mapping  associating to each polynomial $P$ the 
number $deg(P)+1$ if $P$ is not null, and $0$ else. Then $q$ is an ultrametric semi-norm on the 
vector space $\IK[X]$ over the valued field $(\IK,|.|_{disc})$. Now  the ``identity 
transformation'' $\Id: (\IK[X],p) \to (\IK[X],q)$ is continuous (because the topology of the 
semi-normed space $(\IK[X],p)$ is discrete), but $\Id$ is not bounded  
with respect to the  semi-norms $p$ and $q$, since for every $n \in \IN$, $p(X^n)=1$ and 
$q(X^n)=n+1$.
\end{remark}

Given a valued field $(\IK,|.|)$, and two semi-normed $\IK$-vector spaces $(E,p)$ and $(F,q)$, 
we denote by $\BL(E,F)$ the vector space of bounded  linear  mappings from $E$ to $F$.  
Given some bounded  linear mapping
$T: E \to F$, the real number $\inf  \{ M \in \IR_+ : \forall x \in E \; q(T(x))\le M p(x)\}$ is 
called the {\em semi-norm} 
of the operator $T$, and is denoted by $\norme{T}_{\BL(E,F)}$ (or $\norme{T}$). The mapping
$\norm : \BL(E,F) \to \IR_+$ associating to each bounded  operator $T \in \BL(E,F)$ its semi-norm 
$\norme{T}$ is a vector semi-norm, which is ultrametric if the semi-norm $q$ of $F$ is ultrametric. 

\begin{remark}
 Given a spherically complete ultrametric valued field $(\IK,|.|)$, the Ingleton statement 
 $\I_{(\IK,|.|)}$ can be reformulated as follows: for every ultrametric semi-normed space $(E,p)$
 over the valued field $(\IK,|.|)$, for every vector subspace $F$ of $E$ and for every bounded
 linear mapping $f: (F,p) \to (\IK,|.|)$, there exists a bounded linear mapping 
 $\tilde f: (E,p) \to (\IK,|.|)$ extending $f$ such that $\norme{\tilde f}=\norme{f}$.
\end{remark}

Given a valued field $(\IK,|.|)$, a semi-normed $\IK$-vector space $(E,p)$ and a vector 
subspace $F$ of $E$, a {\em continuous linear extender} from $\BL(F,\IK)$ to $\BL(E,\IK)$ 
is a continuous linear mapping $T: \BL(F,\IK)  \to \BL(E,\IK)$ such that for every $f \in 
\BL(F,\IK)$, $T(f)$ extends $f$;
moreover, if for every $f \in \BL(F,\IK)$, $T(f)$ has the same semi-norm as $f$, then the 
continuous 
linear extender $T$ is said to be {\em isometric}.

\subsection{Orthogonal basis  of a finite dimensional ultrametric semi-normed space}
\begin{lemma*}[{\cite[Ex.~ 3.R p.~63]{vanRooij}}] \label{lem:orthocompl}
 Let $(\IK,|.|)$ be a  spherically complete ultrametric valued field. Let $E$ be a $\IK$-vector 
space endowed with a semi-norm $p: E \to \IR$. Given two vector  subspaces $F$ and $G$ of 
$E$ such that $F \oplus G=E$, and denoting by $P_F: E \to F$ and $P_G: E \to G$ the associated 
projections, the following statements are equivalent:
\begin{enumerate}
 \item For every $x \in E$ $p(P_F(x)) \le p(x)$ ({\em i.e.} $P_F$ is bounded and $\norme{P_F} \le 
1$)
\item For every $x \in E$ $p(P_G(x)) \le p(x)$ ({\em i.e.} $P_G$ is bounded and $\norme{P_G} \le 
1$)
\item For every $x \in F$ and every $x \in G$, $p(x_F \oplus x_G)=\max(p(x_F),p(x_G))$.
\end{enumerate}
\end{lemma*}

\begin{definition}
 Given an ultrametric valued field $(\IK,|.|)$, a $\IK$-vector space $E$ and an ultrametric 
semi-norm 
$p: E \to \IR$, two vector subspaces $F$ and $G$ of $E$ satisfying the conditions of 
Lemma~\ref{lem:orthocompl} are said to be {\em orthocomplemented}.
\end{definition}

\begin{lemma*}
 Let $(\IK,|.|)$ be a  spherically complete ultrametric valued field. Let $E$ be a 
finite dimensional $\IK$-vector 
space endowed with an ultrametric semi-norm $p: E \to \IR$. Every one-dimensional vector subspace 
$D$ of $E$ 
has an orthocomplemented subspace in $E$.
\end{lemma*}
\Pr Let $D$ be a one-dimensional vector subspace of $E$. If   $p_{\restriction D}$ is null, then 
every vector subspace $H$ of $E$ such that $H \oplus D=E$ is an orthocomplement of $D$ in $E$. 
Assume that $p_{\restriction D}$ is not null. Let $a \in D \backslash \{0\}$. Let $f: D \to \IK$ 
be the linear mapping such that $f(a)=1$: then for every $x \in D$, 
$|f(x)| \le \frac{1}{p(a)} p(x)$. Since 
the vector space $E$ is finite dimensional, using Remark~\ref{rem:ing_basis},  there exists in 
$\ZF$ a 
linear 
mapping $\tilde f: E \to \IK$ extending $f$ such that $|\tilde f| \le \frac{1}{p(a)} p$. 
Then $H:=\ker(\tilde f)$ is an orthocomplement of the subspace $D$ in $E$: for every $x \in E$,
$x=(x-\tilde f(x).a) \oplus \tilde f(x).a$ where $x-\tilde f(x).a \in H$; moreover, 
$p(\tilde f(x).a)=|\tilde f(x)| p(a) \le p(x)$, so, denoting by $P_D$ the projection onto the 
subspace  $D$ with kernel $H$, $\norme{P_D} \le 1$.
\Rp 

\begin{remark}
 Given a spherically complete ultrametric valued field $(\IK,|.|)$, and an ultrametric semi-normed 
$\IK$-vector space $(E,p)$, then the Ingleton statement $\I_{(\IK,|.|)}$ implies that every 
one-dimensional subspace of $E$ is orthocomplemented in $E$.
\end{remark}

\begin{definition}
 Given an ultrametric valued field $(\IK,|.|)$, a $\IK$-vector space $E$ and an ultrametric 
semi-norm 
$p: E \to \IR$, a  sequence $(e_i)_{0 \le i \le p}$ of $E$ is said to be $p$-orthogonal if for 
every 
sequence $(s_i)_{0 \le i \le p}$ of $\IK$, 
$p(\sum_{0 \le i \le p}s_i.e_i)=\sup_{0 \le i \le p}p(s_i.e_i)$.
\end{definition}

\begin{lemma*}[{\cite[Lemma~ 5.3 p.169]{vanRooij}}] 
 Let $(\IK,|.|)$ be a  spherically complete ultrametric valued field. For every finite 
dimensional $\IK$-vector space $E$ and every ultrametric semi-norm $p: E \to \IR$, there exists a 
$p$-orthogonal basis in the vector space $E$. 
\end{lemma*}
\Pr  We prove the Lemma by recursion over the dimension of $E$. If $dim(E)=1$, then every basis of 
$E$ is $p$-orthogonal.  Assume that  the result holds for some natural number $n \ge 1$ and 
assume that $E$ is a $\IK$-vector space with dimension $n+1$, and that $p: E \to \IR$ is a 
ultrametric semi-norm. Let $a \in E \backslash \{0\}$, and let $D$ be the line $\IK.a$. Using the 
previous Lemma, let $H$ be a $p$-orthocomplemented subspace of $D$ in $E$. Using the recursion 
hypothesis, let $(e_i)_{1 \le i \le n}$ be a $p$-orthogonal basis of the subspace $H$. Let 
$e_{n+1}:=a$. Since $D$ and $H$ are orthocomplemented, $(e_i)_{0 \le i \le n+1}$ is a 
$p$-orthogonal basis of $E$. 
\Rp

\subsection{Ultrametric isometric linear extenders}
In this Section, we shall show that given a spherically complete valued field $(\IK,|.|)$, the 
statement $\I_{(\IK,|.|)}$ is equivalent to the following  ``isometric linear 
extender'' statement:\\
$\LE_{(\IK,|.|)}$: {\em ``For every vector subspace $F$ of an ultrametric semi-normed 
$\IK$-vector space $(E,p)$, there exists an isometric linear extender $T:\BL(F,\IK) \to 
\BL(E,\IK)$.}
\begin{theorem}
Given a spherically complete ultrametric valued field $(\IK,|.|)$, the statements
 $\I_{(\IK,|.|)}$ and $\LE_{(\IK,|.|)}$ are equivalent. 
\end{theorem}
\Pr $\I_{(\IK,|.|)} \Rightarrow \LE_{(\IK,|.|)}$ Let $(E,p)$ be an  ultrametric semi-normed vector 
space $(E,p)$ over the valued 
field $(\IK,|.|)$. We endow the $\IK$-vector space 
$\BL(E,\IK)$ with its ultrametric semi-norm $\norm$. 
Given some vector subspace $F$ of $E$, let $I$ be the set of mappings 
$\Phi: \BL(F,\IK) \to \BL(E,\IK)$ and let $R$ be the binary relation on $fin(\BL(F,\IK)) \times I$ 
such that 
for every $Z  \in fin(\BL(F,\IK))$ and every $\Phi \in I$, $R(Z,\Phi)$ if and only if for every
$f \in Z$, the bounded  linear form  $\Phi(f)$ extends $f$, $\norme{\Phi(f)}=\norme{f}$, 
and $\Phi$ is $\IK$-linear on $\spanv_{\BL(F,\IK)}(Z)$.
Then the binary relation $R$ is concurrent: given $m$ finite subsets 
$Z_1, \dots, Z_m  \in fin(\BL(F,\IK))$, let $B=\{f_1,\dots,f_n\}$ be a (finite) $\norm$-orthogonal 
basis of the $\IK$-vector 
subspace of $\BL(F,\IK)$ generated by the finite set $\cup_{1 \le i \le m}Z_i$; 
then using  $\I_{(\IK,|.|)}$,  let 
$\tilde f_1$, \dots, $\tilde f_n$ be bounded   linear forms on $E$ extending 
$f_1$, \dots, $f_n$ such that for each $i$, $\norme{\tilde f_i}=\norme{f_i}$, and
let $L: \spanv(\{f_1,\dots,f_n\}) \to \BL(E,\IK)$ be the linear mapping such that for each 
$i \in \{1,\dots,n\}$, $L(f_i)=\tilde{f_i}$. For every $f \in \spanv(\{f_1,\dots,f_n\})$,
let us check that $\norme{L(f)}=\norme{f}$. Given $f \in \spanv(\{f_1,\dots,f_n\})$, 
$f$ is of the form $\sum_{1 \le i \le n}s_if_i$ where  $s_1$, \dots, $s_n \in \IK$,
thus $L(f)=\sum_{1\le i \le n}s_i \tilde{f_i}$; since $\norm$ is ultrametric, it follows that  
$\norme{L(f)}=
\norme{\sum_{1\le i \le p}s_i \tilde{f_i}} \le \max_{1 \le i \le n}  \norme{s_i \tilde{f_i}}
=\max_{1 \le i \le n}  |s_i| \norme{ \tilde{f_i}}
=\max_{1 \le i \le n} |s_i|\norme{f_i}
=\max_{1 \le i \le n} \norme{s_i f_i}=\norme{\sum_{1\le i \le n}s_i f_i}$ (because the sequence 
$(f_i)_{1 \le i \le n}$ is $\norm$-orthogonal).
Let $\Phi: \BL(F,\IK) \to \BL(E,\IK)$ be some mapping extending $L$ (for example, define $\Phi(f)=0$
for every $f \in \BL(F,\IK) \backslash \spanv(\{f_1,\dots,f_n\})$. 
Then for every $i \in \{1, \dots,m\}$, $R(Z_i,\Phi)$. 
Consider the filter $\mathcal F$ on $I$ generated by the set $\{R(Z) ; Z \in fin(\BL(F,\IK))\}$. 
Then the mapping  $\Phi : \BL(F,\IK) \to L(E,\IK_{\mathcal F})$
 associating to each $f \in \BL(F,\IK)$ the $\IK$-linear mapping 
$\Phi(f):  E \to {\IK}_{\mathcal F}$  associating to each  
$x \in E$ the class of  $(f(x))_{f \in  I}$ in  $ {\IK}_{\mathcal F}$
 is linear, and for every $f \in \BL(F,\IK)$, $\Phi(f) :E \to {\IK}_{\mathcal F}$ extends $f$ and 
for every $x \in E$, $N_{|.|,\mathcal F}(\Phi(f)(x)) \le \norme{f} \norme{x}$. Using 
$\I_{(\IK,|.|)}$ and Lemma~\ref{lem:equiv-ing}, consider   a  $\IK$-linear retraction 
$r: ({\IK}_{\mathcal F})_b \to \IK$ of 
$j_{\mathcal F}: \IK \to {\IK}_{\mathcal F}$ such that for every 
$x \in ({\IK}_{\mathcal F})_b$, $|r(x)| \le N_{|.|,\mathcal F}(x)$. 
Thus, for every $f \in \BL(F,\IK)$, for every $x \in F$, 
$|r \circ {\Phi(f)}(x)| \le \norme{f} \norme{x}$,
so $\norme{r \circ {\Phi(f)}} \le \norme{f}$; since $r \circ {\Phi(f)}: E \to \IK$ extends $f: F 
\to \IK$, it follows that $\norme{r \circ {\Phi(f)}} = \norme{f}$. 
Let $T: \BL(F,\IK) \to \BL(E,\IK)$ be the mapping 
$f \mapsto r \circ {\Phi(f)}$. 
Then the mapping $T: \BL(F,\IK) \to \BL(E,\IK)$ is an isometric linear extender.\\
The implication  $\LE_{(\IK,|.|)} \Rightarrow  \I_{(\IK,|.|)}$ is trivial.
\Rp

\bibliographystyle{abbrv}
\bibliography{linear_extenders_rvst1.bbl}
\end{document}